\documentclass{amsart}
\usepackage{amssymb}

\newtheorem{theorem}{Theorem}

\newcommand{\C}{\mathbb{C}}

\newcommand{\Z}{\mathbb{Z}}
\newcommand{\N}{\mathbb{N}}
\newcommand{\PP}{\mathbb{P}}
\newcommand{\LL}{\mathcal{L}}
\newcommand{\CC}{\mathcal{C}}
\newcommand{\RR}{\mathcal{R}}
\newcommand{\LLL}{\mathbf{L}}
\newcommand{\OO}{\mathcal{O}}
\newcommand{\ox}{\vec{x}}

\newcommand{\oc}{\vec{c}}
\newcommand{\XX}{\mathbb{X}}

\newcommand{\gf}{\mathfrak{g}}
\newcommand{\kdelta}{\mbox{\boldmath $\delta$}}
\newcommand{\dbslash}{/\!\!/}

\DeclareMathOperator{\ind}{ind}
\DeclareMathOperator{\ad}{ad}
\DeclareMathOperator{\GL}{GL}
\DeclareMathOperator{\gl}{gl}
\DeclareMathOperator{\coh}{coh}

\DeclareMathOperator{\Rep}{Rep}
\DeclareMathOperator{\Mat}{Mat}

\DeclareMathOperator{\rank}{rank}
\DeclareMathOperator{\iso}{iso}
\DeclareMathOperator{\modu}{mod}
\DeclareMathOperator{\dimv}{\underline{\dim}}

\title[Quiver algebras]{Quiver algebras, weighted projective lines, and the Deligne-Simpson problem}
\author{William Crawley-Boevey}
\address{Department of Pure Mathematics, University of Leeds, Leeds LS2 9JT, UK}\email{w.crawley-boevey@leeds.ac.uk}

\begin{document}
\begin{abstract}
We describe recent work on preprojective algebras and moduli
spaces of their representations.
We give an analogue of Kac's Theorem,
characterizing the dimension types of indecomposable
coherent sheaves over weighted projective lines in terms of
loop algebras of Kac-Moody Lie algebras, and explain how
it is proved using Hall algebras.
We discuss applications to the problem of describing the
possible conjugacy classes of sums and products of matrices
in known conjugacy classes.
\end{abstract}



\maketitle
\section*{Introduction}
\label{s:intro}
Preprojective algebras were introduced by Gelfand and Ponomarev,
and in a deformed version by Crawley-Boevey and Holland.
They arose in the theory of representations
of quivers, but have interesting links with
Kleinian singularities, Kac-Moody Lie algebras and
noncommutative symplectic geometry.
In the first part of this article, \S\ref{s:preproj},
we survey some of the results we have obtained in the
last ten years concerning these algebras,
and moduli spaces of their representations.

The Deligne-Simpson problem asks about the existence of matrices in given
conjugacy classes, with product the identity, and no common invariant subspace.
Some time ago it became clear that our work on preprojective
algebras solves an additive analogue.
To solve the original problem, one needs to
pass to a new setup, in which representations
of quivers are replaced by coherent sheaves on
weighted projective lines (or parabolic bundles),
and representations of the
preprojective algebra are replaced by logarithmic
connections. We discuss all this in~\S\ref{s:dsp}.

A key ingredient in the theory of preprojective algebras
is Kac's Theorem, describing the possible dimension
vectors of indecomposable representations of quivers.
In the new setup, one needs an analogue of Kac's Theorem
for weighted projective lines.
We discuss it in \S\ref{s:wtd}, and
outline a proof via Hall algebras in~\S\ref{s:hall}.

In the rest of this introduction we recall some basic
facts about representations of quivers.
A \emph{quiver} $Q$, or more precisely $(I,Q,h,t)$, consists of finite sets
$I$ and $Q$ of vertices and arrows, and maps $h,t:Q\to I$,
assigning to each arrow its head and tail vertices.
We fix a base field $K$, algebraically closed unless otherwise indicated.
By a \emph{representation} $X$ of $Q$, one means the assignment of a vector space
$X_v$ for each vertex $v$, and a linear map $X_{t(a)}\to X_{h(a)}$ for each arrow~$a$.
There are natural notions of homomorphisms between representations,
sub-representations, etc.

The \emph{path algebra} $KQ$ has basis
the paths $a_1 a_2\dots a_n$ in $Q$ of length $n\ge 1$,
with $t(a_i) = h(a_{i+1})$ for all $i$,
and a \emph{trivial path} of length 0 for each vertex~$v$.
It is an associative algebra, with the
product of two paths given by their concatenation,
if this makes sense, and otherwise zero, and
the sum of the trivial paths is a multiplicative identity.
The category of representations of $Q$ is equivalent
to the category of left $KQ$-modules, so one can
use homological algebra, composition series,
the Krull-Remak-Schmidt Theorem, and so on.

We now assume for simplicity that $Q$ has no oriented cycles,
in which case $KQ$ is finite dimensional,
although many results hold without this restriction.

Let $\gf$ be the Kac-Moody Lie algebra
given by the symmetric generalized Cartan
matrix $(a_{uv})_{u,v\in I}$ which has diagonal entries 2,
and off-diagonal entries given by minus the number
of arrows in $Q$ between $u$ and $v$, in either direction.
Thus $\gf$ is generated over $\C$ by $e_v,f_v,h_v$ ($v\in I$)
with relations
\begin{equation}
\begin{cases}
\label{e:kacrelns}
[h_u,h_v] =0,
\quad
[e_u,f_v]=\kdelta_{uv} h_v,
\\
[h_u,e_v] = a_{uv} e_v,
\quad
[h_u,f_v] = -a_{uv} f_v,
\\
(\ad e_u)^{1-a_{uv}}(e_v)=0,
\quad
(\ad f_u)^{1-a_{uv}}(f_v)=0
\quad
(\text{if $u\neq v$}),
\end{cases}
\end{equation}
where $\kdelta$ is the Kronecker delta function.
The root lattice $\Gamma$ of $\gf$ is the free additive group on
symbols $\alpha_v$ ($v\in I$), it grades $\gf$, with
$\deg e_v = \alpha_v$, $\deg f_v=-\alpha_v$ and $\deg h_v=0$,
and the set of roots is
$\Delta = \{ 0\neq \alpha\in \Gamma \mid \gf_\alpha\neq 0 \}$.
Recall that there are real roots, obtained from the simple
roots $\alpha_u$ by a sequence of reflections
$s_v(\alpha) = \alpha-(\alpha,\alpha_v)\alpha_v$,
where $(-,-)$ is the symmetric bilinear
form on $\Gamma$ with $(\alpha_u,\alpha_v) = a_{uv}$,
and there may also be imaginary roots.
Defining $p(\alpha) = 1 - \frac12 (\alpha,\alpha)$,
the real roots have $p(\alpha)=0$, and the imaginary roots have
$p(\alpha)>0$.

Gabriel's Theorem \cite{Gabriel} asserts that a quiver $Q$ has only finitely
many indecomposable representations if and only if $\gf$ is of finite type,
i.e.\ the underlying graph of $Q$ is a Dynkin diagram (of type ADE).
In this case the map sending a representation to its dimension vector
\[
\dimv X = \sum_{v\in I} (\dim X_v) \alpha_v \in \Gamma
\]
gives a 1-1 correspondence between  indecomposable representations
and positive roots.
Kac's Theorem \cite{Kac1,Kac2} extends this to $\gf$ of arbitrary type:
the dimension vectors of indecomposable representations
are exactly the positive roots, there is a unique
indecomposable for each real root, infinitely many for each imaginary root.

\section{Preprojective algebras}
\label{s:preproj}
The \emph{preprojective algebra} of a quiver $Q$ is the algebra
\[
\Pi(Q) = K \overline{Q} /
\biggl( \sum_{a\in Q} (aa^* - a^*a) \biggr),
\]
where the \emph{double} $\overline{Q}$ of $Q$ is obtained
by adjoining a reverse arrow $a^*$ for each arrow $a\in Q$.
In the finite type case it is isomorphic, as a $KQ$-module,
to the direct sum of one copy of each indecomposable representation of $Q$.
In general, it is isomorphic to the direct sum of the
indecomposable representations of $Q$ that are \emph{preprojective},
meaning that some power of the Coxeter functor, or equivalently of the
Auslander-Reiten translation $DTr$, sends them to a projective $KQ$-module.

Preprojective algebras first appeared
in unpublished work of I.~M.~Gelfand and V.~A.~Ponomarev,
in a lecture delivered by A.~V.~Roiter at the
Second International Conference on
Representations of Algebras (Ottawa, 1979).
See \cite{Ringelpreproj} for a discussion
about variations on this definition.
Note also that work by Riedtmann \cite{Riedtmann} contains 
parallel ideas.

Given $\lambda = (\lambda_v)_{v\in I} \in K^I$,
Crawley-Boevey and Holland \cite{CBH} have introduced
the \emph{deformed preprojective algebra}, $\Pi^\lambda(Q)$,
in which the relation is replaced by
\begin{equation}
\label{e:defrel}
\sum_{a\in Q} (aa^* - a^*a) - \lambda,
\end{equation}
where $\lambda$ is identified with the corresponding
linear combination of trivial paths.

Up to isomorphism, these algebras do not depend on the orientation of~$Q$.
They are related to some elementary symplectic geometry.
Choosing bases for the vector spaces, representations
of $Q$ of dimension vector $\alpha = \sum_{v} n_v \alpha_v$
are given by elements of the space
\[
\Rep(Q,\alpha) = \bigoplus_{a\in Q} \Mat_{n_{h(a)}\times n_{t(a)}}(K),
\]
and isomorphism classes correspond to orbits of the group
\[
\GL(\alpha) = \prod_{v\in I} \GL_{n_v}(K)
\]
acting by conjugation.
The space of representations of $\overline Q$ can then be identified
with a cotangent bundle
\[
\Rep(\overline{Q},\alpha) \cong \Rep(Q,\alpha) \times \Rep(Q,\alpha)^* \cong T^*\Rep(Q,\alpha).
\]
This has a natural symplectic structure, and
associated to the action of $\GL(\alpha)$ there is a moment map
\[
\mu_\alpha : \Rep(\overline{Q},\alpha) \to \gl(\alpha),
\quad
x \mapsto \biggl(\sum_{\substack{a\in Q\\ h(a)=v}} x_a x_{a^*}
- \sum_{\substack{a\in Q\\ t(a)=v}} x_{a^*} x_a \biggr)_{v\in I}\ .
\]
Identifying $\lambda$ with a central element of $\gl(\alpha)$,
there is a quotient
\[
N_Q(\lambda,\alpha) = \mu_\alpha^{-1}(\lambda)\dbslash \GL(\alpha),
\]
a `symplectic quotient', or `Marsden-Weinstein reduction'.
Here the double slash denotes the affine variety
that classifies closed orbits of $\GL(\alpha)$ on $\mu_\alpha^{-1}(\lambda)$.
Now the elements of this fibre are exactly
the representations of $\overline{Q}$ satisfying (\ref{e:defrel}),
so $N_Q(\lambda,\alpha)$ classifies isomorphism classes of
semisimple representations of $\Pi^\lambda(Q)$ of dimension vector~$\alpha$.

Since the trace of $\mu_\alpha(x)$ is always zero,
if there is a representation of $\Pi^\lambda(Q)$
of dimension vector $\alpha = \sum n_v \alpha_v$, then
$\lambda\cdot\alpha = \sum \lambda_v n_v$ must be zero.
If $Q$ is a Dynkin diagram, then $\Pi^\lambda(Q)$ is finite-dimensional,
and this argument shows that for generic $\lambda$ it is even the zero algebra.

The case when $Q$ is an extended Dynkin diagram, or
equivalently when $\gf$ is an affine Lie algebra, appeared
in work of Kronheimer~\cite{Kronheimer},
made more explicit by Cassens and Slodowy~\cite{CS}.
If $\delta$ is the minimal positive imaginary root,
then $N_Q(0,\delta)$ is the corresponding Kleinian surface singularity,
and the spaces $N_Q(\lambda,\delta)$, for suitably varying $\lambda$, give its
semiuniversal deformation.
The key idea of \cite{CBH} is that
the deformed preprojective algebras $\Pi^\lambda(Q)$, for unrestricted
$\lambda$, give a larger family of deformations of the Kleinian
singularity, the general one being noncommutative.

\begin{theorem}[Crawley-Boevey and Holland \cite{CBH}]
Suppose $Q$ is an extended Dynkin diagram,
and $e$ is the trivial path corresponding to an
extending vertex for $Q$.
Then $\OO^\lambda = e \Pi^\lambda(Q) e$
is a noetherian domain of Gelfand-Kirillov dimension~2,
Auslander-Gorenstein and Cohen-Macaulay.
Moreover, it is commutative if and only if $\lambda\cdot\delta=0$,
and if so, it is isomorphic to the coordinate ring of $N_Q(\lambda,\delta)$.
\end{theorem}

For some further work on the $\OO^\lambda$, see \cite{BGK,Boyarchenko,Holland}.
Returning to the general case, to decide when $N_Q(\lambda,\alpha)$ is nonempty,
one needs to know whether or not
there is a representation of $\Pi^\lambda(Q)$ of dimension
vector $\alpha$.
It is not hard to show that
a representation $X$ of $Q$ is in the image of the projection
$\mu_\alpha^{-1}(\lambda)\to \Rep(Q,\alpha)$,
so is the restriction of a representation
of $\Pi^\lambda(Q)$, if and only if the
dimension vector $\beta$ of each
indecomposable direct summand of $X$ satisfies $\lambda\cdot\beta=0$.
With Kac's Theorem, this gives the following.

\begin{theorem}[\cite{CBmm}]
\label{t:pilamdim}
The space $N_Q(\lambda,\alpha)$ is nonempty, or equivalently
there is a representation of $\Pi^\lambda(Q)$ of dimension vector $\alpha$,
if and only if $\alpha$ can be written as a sum of positive roots
$\alpha = \beta+\gamma+\dots$ with
$\lambda\cdot\beta = \lambda\cdot\gamma = \dots = 0$.
\end{theorem}

More difficult is the following.

\begin{theorem}[\cite{CBmm}]
\label{t:pilamsim}
There is a simple representation of $\Pi^\lambda(Q)$ of dimension vector $\alpha$
if and only if $\alpha$ is a positive root, $\lambda\cdot\alpha=0$, and
$p(\alpha) > p(\beta)+p(\gamma)+\dots$ for any nontrivial decomposition
of $\alpha$ as a sum of positive roots $\alpha = \beta+\gamma+\dots$
with $\lambda\cdot\beta = \lambda\cdot\gamma = \dots = 0$.
\end{theorem}

We have also shown \cite{CBdecomp,CBnorm} that
if $K$ has characteristic zero and $N_Q(\lambda,\alpha)$ is nonempty,
then it is an irreducible normal variety.
Bocklandt and Le Bruyn \cite{BLB,LeBruyn} have obtained further
results in this direction.
See \cite{CBEG} for more about the link with noncommutative symplectic geometry.

There are more general moduli spaces $N_Q(\lambda,\alpha)_\theta$,
depending on suitable stability data~$\theta\in \Z^I$.
Examples of these are the `quiver varieties' used by Nakajima to
construct integrable representations of Lie algebras and quantum
groups, see his ICM talk~\cite{Nakajimaicm} (and \cite{CBmm}).
Note that over $\C$, the moduli spaces are special cases
of hyper-K\"ahler quotients, and by a standard trick of changing
the complex structure, $N_Q(0,\alpha)_\theta$
is homeomorphic to $N_Q(\theta,\alpha)$.

We now give an application of these ideas.
When working over a finite field, it is natural to consider
representations of $Q$ which are \emph{absolutely indecomposable},
meaning that they remain indecomposable over the algebraic closure
of the field.
Kac showed that up to isomorphism,
the number such representations of dimension vector $\alpha$
is polynomial in the size $q$ of the field,
of the form $a_\alpha(q)$ for some $a_\alpha\in\Z[t]$.
He conjectured that $a_\alpha$ has non-negative coefficients,
and that the constant term is the root multiplicity $\dim \gf_\alpha$.
In partial answer we have the following.

\begin{theorem}[Crawley-Boevey and Van den Bergh \cite{CBVdB}]
\label{t:kacconj}
If $\alpha = \sum_{v\in I} n_v \alpha_v$ is indivisible,
meaning that the $n_v$ have no common divisor, then
$a_\alpha$ has non-negative coefficients, and the constant
term is the root multiplicity $\dim \gf_\alpha$.
\end{theorem}

We explain the positivity.
Since $\alpha$ is indivisible,
one can fix $\lambda\in\Z^I$ with $\lambda\cdot\alpha=0$
but $\lambda\cdot\beta\neq 0$ for all $0<\beta<\alpha$.
The argument of Theorem~\ref{t:pilamdim}
shows that the number of points in $N_Q(\lambda,\alpha)$, over
a field with $q$ elements and sufficiently large characteristic, is $q^{p(\alpha)} a_\alpha(q)$,
and then if $N_Q(\lambda,\alpha)$ had been a projective variety,
the Weil conjectures would have given positivity.
However, $N_Q(0,\alpha)_\lambda$ is sufficiently close to being
projective for the Weil conjectures to apply to it, and
by the hyper-K\"ahler trick, the cohomologies of
$N_Q(0,\alpha)_\lambda$ and $N_Q(\lambda,\alpha)$
are isomorphic when the base field is the algebraic closure of a finite field
of sufficiently large characteristic. Moreover, it is possible
to ensure that this isomorphism is compatible with Frobenius
maps, so that $N_Q(\lambda,\alpha)$ is good enough.

\section{Weighted projective lines}
\label{s:wtd}
In this section we give an analogue of Kac's Theorem
for weighted projective lines. When studying representations of finite-dimensional
associative algebras, quivers tell one about \emph{hereditary} algebras,
i.e.\ those with global dimension $\le 1$.
One of the breakthroughs in this area
was the discovery by Brenner and Butler \cite{BB,HR}
of algebras $A$ that are `tilted' from a hereditary algebra $H$.
In Happel's language~\cite{Happelbook},
there is a derived equivalence
$D^b(\modu A) \simeq D^b(\modu H)$,
which is useful since in $D^b(\modu H)$
any indecomposable object is a shift of a module.

Geigle and Lenzing \cite{GL} realized that
there are other algebras,
including Ringel's `canonical algebras' \cite{Ringelbook},
which aren't necessarily tilted from
hereditary algebras, but are tilted from suitable
hereditary abelian categories.
See Reiten's ICM talk \cite{Reitenicm}
for further progress in this direction.

We concentrate on Geigle and Lenzing's categories.
A \emph{weighted projective line} $\XX$ is specified
by giving a collection of distinct points $D=(a_1,\dots,a_k)$ in
the projective line $\PP^1$ over $K$, and a
\emph{weight sequence} $\mathbf{w}=(w_1,\dots,w_k)$,
that is, a sequence of positive integers.
The category, $\coh\XX$, of coherent sheaves on $\XX$,
can be defined as the quotient of the category of finitely
generated $\LLL(\mathbf{w})_+$-graded $S(\mathbf{w},D)$-modules
by the Serre subcategory of finite length modules. Here
$\LLL(\mathbf{w})$ is the additive group with
generators $\ox_1,\dots,\ox_k,\oc$ and relations $w_1 \ox_1 = \dots = w_k \ox_k = \oc$,
partially ordered, with positive cone
$\LLL(\mathbf{w})_+ =  \N \oc+\sum_{i=1}^k \N \ox_i$,
and
\[
S(\mathbf{w},D) = K[u,v,x_1,\dots,x_k] / ( x_i^{w_i} - \lambda_i u - \mu_i v ),
\]
where $a_i = [\lambda_i:\mu_i]\in\PP^1$,
with grading $\deg u = \deg v = \oc$ and $\deg x_i = \ox_i$.
Geigle and Lenzing showed that $\coh\XX$ is a hereditary abelian category;
the free module gives a structure sheaf $\OO$, and shifting the grading gives
twists $E(\ox)$ for any sheaf~$E$ and $\ox\in\LLL(\mathbf{w})$;
also, every sheaf is the direct sum of a `torsion-free' sheaf,
with a filtration by sheaves of the form $\OO(\ox)$,
and a finite-length `torsion' sheaf.

Let $Q_{\mathbf{w}}$ be the star-shaped quiver
whose vertex set $I$ consists of a central
vertex~$*$, and vertices, denoted $ij$ or $i,j$, for
$1\le i\le k$, $1\le j<w_i$, and with arrows
\mbox{$*\leftarrow i1 \leftarrow i2 \leftarrow \dots$} for all $i$.
The appropriate Lie algebra to consider is the
loop algebra $L\gf = \gf [t,t^{-1}]$,
where $\gf$ is the Kac-Moody algebra associated $Q_{\mathbf{w}}$,
or, better, an extension $\LL\gf$ with generators
$e_{vr},f_{vr},h_{vr}$ ($v\in I$, $r\in \Z$) and $c$ subject to the relations
\begin{equation}
\label{e:loopalg}
\begin{cases}
\text{$c$ central,}
\quad
[e_{vr},e_{vs}] = 0,
\quad
[f_{vr},f_{vs}] = 0,
\\
[h_{ur},h_{vs}]=ra_{uv}\,\kdelta_{r+s,0}\,c,
\quad
[e_{ur}, f_{vs}] = \kdelta_{uv}\, \left(h_{v,r+s} + r\,\kdelta_{r+s,0}\,c\right),
\\
[h_{ur},e_{vs}] = a_{uv} e_{v,r+s},
\quad
[h_{ur},f_{vs}] = -a_{uv} f_{v,r+s},
\\
(\ad e_{u0})^{1-a_{uv}}(e_{vs}) = 0,
\quad
(\ad f_{u0})^{1-a_{uv}}(f_{vs}) = 0
\quad
(\text{if $u\neq v$}),
\end{cases}
\end{equation}
see~\cite{MRY}. The root lattice for either algebra is
$\hat{\Gamma}=\Gamma\oplus\Z \delta$, with
$\deg e_v t^r = \deg e_{vr} = \alpha_v + r\delta$,
$\deg f_v t^r = \deg f_{vr} = -\alpha_v + r\delta$,
$\deg h_v t^r = \deg h_{vr} = r\delta$
and $\deg c=0$,
and the set of roots for either algebra is
\[
\hat\Delta = \{ \alpha+r\delta \mid \alpha\in\Delta, r\in\Z \} \cup \{ r\delta \mid 0\neq r\in Z \}.
\]
The real roots are $\alpha+r\delta$ with $\alpha$ real.
(If $\gf$ is of finite type, $\LL\gf$ is the corresponding
affine Lie algebra, and if $\gf$ is of affine type,
$\LL\gf$ is a toroidal algebra.)

The Grothendieck group $K_0(\coh\XX)$ was computed by Geigle and Lenzing,
and following Schiffmann~\cite{Schiffmann} it can be
identified with $\hat{\Gamma}$.
Now $K_0(\coh\XX)$ is partially ordered, with the positive cone being the classes
of objects in $\coh\XX$, and this gives a partial ordering on $\hat \Gamma$.

\begin{theorem}[\cite{CBkacw}]
\label{t:kaccoh}
If $\XX$ is a weighted projective line,
there is an indecomposable coherent sheaf on $\XX$ of type $\phi\in\hat \Gamma$ if
and only if $\phi$ is a positive root.
There is a unique indecomposable for a real root, infinitely many for an imaginary root.
\end{theorem}

We remark that there is a classification of the indecomposables
if $\gf$ is of finite type \cite{GL},
or affine type \cite{LM}. The latter
is essentially equivalent to the classification
for tubular algebras,
see Ringel's ICM talk \cite{Ringelicm,Ringelbook}.

Lenzing \cite[\S 4.2]{Lenzing} showed that the
category of torsion-free sheaves on $\XX$ is equivalent to the category
of (quasi) parabolic bundles on $\PP^1$ of weight type $(D,\mathbf{w})$, that is,
vector bundles $\pi:E\to \PP^1$ equipped with a flag of vector subspaces
\[
\pi^{-1}(a_i) \supseteq E_{i1} \supseteq \dots \supseteq E_{i,w_i-1}
\]
for each $i$.
This equivalence is not unique, but it can be chosen so that
if $E$ is a parabolic bundle, then $[E] = \dimv E + (\deg E)\delta$,
where the dimension vector is
\[
\dimv E = n_* \alpha_*  + \sum_{i=1}^k \sum_{j=1}^{w_i-1} n_{ij} \alpha_{ij} \in \Gamma,
\]
with $n_* = \rank E$ and $n_{ij} = \dim E_{ij}$.
This is necessarily \emph{strict}, meaning
that $n_* \ge n_{i1} \ge n_{i2} \ge \dots \ge n_{i,w_i-1}\ge 0$.
We can now restate Theorem~\ref{t:kaccoh} as follows (see \cite{CBipb}).
For each $d\in\Z$ there is an indecomposable parabolic bundle of
dimension vector $\alpha$ and degree $d$ if and only if $\alpha$
is a strict root for $\gf$. There is a unique indecomposable for a real root,
infinitely many for an imaginary root.

\section{Hall algebras}
\label{s:hall}
In this section we explain the proof of Theorem~\ref{t:kaccoh}.
Let $\CC$ be an abelian category that is \emph{finitary},
meaning that its Hom and Ext spaces are finite sets.
The \emph{Hall algebra} of $\CC$, over a commutative ring $\Lambda$, is
the free $\Lambda$-module
\[
H_\Lambda(\CC) = \bigoplus_{Z\in \iso\CC} \Lambda u_Z,
\]
with basis the symbols $u_Z$, where $Z$ runs through
$\iso\CC$, a set of representatives of the isomorphism classes of $\CC$.
It is an associative algebra with product
\[
u_X u_Y = \sum_{Z\in \iso\CC} F^Z_{XY} u_Z,
\]
where $F^Z_{XY}$ is the number of subobjects $Z'$ of $Z$ with $Z'\cong Y$ and $Z/Z'\cong X$.
In case $\CC$ is the category of finite abelian groups,
or finite abelian $p$-groups, this notion is due to
Steinitz~\cite{Steinitz} and Hall~\cite{Hall}.
The current interest in Hall algebras stems from
the discovery by Ringel \cite{Ringelhallqg} of a relationship between
quantum groups and Hall algebras for categories
of representations of finite-dimensional hereditary algebras
over finite fields.
This quickly influenced the development
of canonical bases, see Lusztig's ICM talk \cite{Lusztigicm}.

How to recover the underlying Lie algebra? Ringel \cite{Ringelhall}
realized, for finite type hereditary algebras over a finite field $K$,
that if one uses a ring $\Lambda$ in which $|K|=1$,
and $\ind \CC$ is the set of indecomposables in $\iso\CC$,
then the $u_Z$ with $Z\in \ind \CC$
generate a Lie subalgebra of $H_\Lambda(\CC)$ with bracket
\[
[u_X, u_Y] = \sum_Z (F^Z_{XY} - F^Z_{YX}) u_Z.
\]
To get at something resembling a semisimple Lie algebra,
not just its positive part,
there is a construction of Peng and Xiao \cite{Hubery,PX},
in which one starts not with an abelian category, but
with a triangulated $K$-category that is \emph{2-periodic},
meaning that the shift functor $T$ satisfies $T^2=1$.

Let $\XX$ be a weighted projective line over a finite field $K$,
whose marked points are all defined over $K$.
The category $\coh\XX$ is still defined and well-behaved, and
Schiffmann \cite{Schiffmann} has studied its Hall algebra.
Applying the construction of Peng and Xiao to
the orbit category $\RR_\XX = D^b(\coh\XX)/(T^2)$, one obtains a Lie algebra
with triangular decomposition
\[
L_\Lambda(\RR_\XX) =
\biggl(\bigoplus_{X\in \ind \coh\XX} \Lambda u_X \biggr)
\oplus
(\Lambda \otimes_\Z \hat \Gamma)
\oplus
\biggl(\bigoplus_{X\in \ind \coh\XX} \Lambda u_{TX} \biggr),
\]
where $\Lambda$ is still a commutative ring in which $|K|=1$. We have the following result.

\begin{theorem}[\cite{CBkacw}]
\label{t:lieelsts}
$L_\Lambda(\RR_\XX)$ contains elements
$e_{vr},f_{vr},h_{vr}$ $(v\in I,r\in \Z)$ and $c$ satisfying
the relations \textup{(\ref{e:loopalg})} for $\LL\gf$.
\end{theorem}

The elements are explicitly given: $c = -1\otimes\delta$,
$e_{*,r} = u_{\OO(r\oc)}$, $f_{*,r} = -u_{T\OO(-r\oc)}$, and
the $e_{ij,r}$ and $f_{ij,r}$ are all of the form $u_X$ or $-u_{TX}$ for
suitable indecomposable torsion sheaves $X$.
See also \cite{LP}, where elliptic Lie algebra generators
are found in $L_\Lambda(\RR_\XX)$ for $\gf$ of affine type.

Concerning the proof of Theorem~\ref{t:kaccoh}, the main problem
is to show that the number of indecomposables of type $\phi=\alpha+r\delta$
is the same as the number of type $s_v(\alpha)+r\delta$.
By arguments already used in the proof of Kac's Theorem,
one reduces to counting numbers of indecomposables
for weighted projective lines over finite fields, so
dimensions of root spaces of $L_\Lambda(\RR_\XX)$.
A standard argument in the theory of complex Lie algebras,
using $\mathrm{sl}_2$-triples $(e,f,h)$, shows that the root multiplicities
for roots related by a reflection are equal.
Now Theorem~\ref{t:lieelsts} provides such triples,
and although the argument uses the fact that the base field has
characteristic zero, for example it involves the operator $\exp(\ad e)$
with $\ad e$ acting locally nilpotently, it works if
$\Lambda$ is a field of sufficiently large characteristic, and this
can be arranged by taking the finite field $K$ to be sufficiently large.

\section{The Deligne-Simpson problem}
\label{s:dsp}
Given invertible matrices in known conjugacy (i.e.\ similarity) classes,
what can one say about the conjugacy class
of their product? More symmetrically,
given conjugacy classes $C_1,\dots,C_k$ in $\GL_n(\C)$,
is there a solution to the equation
\begin{equation}\label{e:prod}
A_1 A_2 \dots A_k = 1
\end{equation}
with $A_i\in C_i$?
The additive analogue asks for a solution to the equation
\begin{equation}
\label{e:sum}
A_1 + A_2 + \dots + A_k = 0,
\end{equation}
where the conjugacy classes may now be in $\gl_n(\C)$.
In full generality these problems seem to be open,
but there are partial results, see for example~\cite{NetoS}.
The former arises when studying linear ODEs
\begin{equation}
\label{e:ode}
\frac{d^n f}{dz^n} + c_1(z) \frac{d^{n-1} f}{dz^{n-1}} + \dots + c_n(z) f = 0
\end{equation}
whose coefficients are rational functions of $z$.
If $D = \{a_1,\dots,a_k\}$ is the set of singular points of the
coefficients in $\PP^1$, the monodromy of (\ref{e:ode})
is a representation in $\GL_n(\C)$
of the fundamental group
of the punctured Riemann sphere $\PP^1\smallsetminus D$, and
the presentation of this group as
$\langle g_1,\dots,g_k \mid g_1 g_2 \dots g_k = 1\rangle$,
where $g_i$ is a suitable loop
around $a_i$, shows how equation (\ref{e:prod}) arises.

To fix the conjugacy classes
we choose a weight sequence $\mathbf{w}=(w_1,\dots,w_k)$,
and a collection of
complex numbers $\xi = (\xi_{ij})$ ($1\le i\le k$, $1\le j\le w_i$) with
$
(A_i - \xi_{i1} 1) (A_i - \xi_{i2} 1) \dots (A_i - \xi_{i,w_i} 1) = 0
$
for $A_i\in C_i$.
Clearly, if one wishes one can take $w_i$ to be
the degree of the minimal polynomial of $A_i$, and
$\xi_{i1},\dots,\xi_{i,w_i}$ to be its roots.
Let $Q_{\mathbf{w}}$ be the quiver associated to
$\mathbf{w}$ as in \S\ref{s:wtd}, let $I$ be its vertex set,
let $\gf$ be the corresponding Kac-Moody Lie algebra,
and let $\Gamma$ be its root lattice.
The $C_i$ determine an element $\alpha = \sum_v n_v \alpha_v\in \Gamma$,
with $n_* = n$ and
\[
n_{ij} = \rank (A_i - \xi_{i1}1)(A_i - \xi_{i2}1)\dots(A_i - \xi_{ij}1)
\]
for $A_i\in C_i$, and conversely $\mathbf{w}$, $\xi$, and $\alpha$ determine the $C_i$.
We define
\[
\xi^{[\beta]} = \prod_{i=1}^k \prod_{j=1}^{w_i} \xi_{ij}^{m_{i,j-1} - m_{ij}},
\quad
\xi*{[\beta]} = \sum_{i=1}^k \sum_{j=1}^{w_i} \xi_{ij}{(m_{i,j-1} - m_{ij})},
\]
for $\beta = \sum m_v \alpha_v$,
with the convention that $m_{i0} = m_*$ and $m_{i,w_i}=0$.
Theorem~\ref{t:pilamdim},
applied to a deformed preprojective algebra $\Pi^\lambda(Q_{\mathbf{w}})$,
gives the following.

\begin{theorem}[\cite{CBipb}]
\label{t:addconjclo}
There is a solution to $A_1+\dots+A_k = 0$ with $A_i$ in the closure $\overline{C_i}$
of $C_i$ if and only if
$\alpha$ can be written as a sum of positive roots $\alpha = \beta+\gamma+\dots$ with
$\xi*{[\beta]} = \xi*{[\gamma]} = \dots = 0$.
\end{theorem}

A solution to equation (\ref{e:prod}) or (\ref{e:sum}) is \emph{irreducible}
if the $A_i$ have no common invariant subspace.
Theorem~\ref{t:pilamsim} gives the following.

\begin{theorem}[\cite{CBad}]
\label{t:ad}
There is an irreducible solution
to $A_1+\dots+ A_k = 0$ with $A_i\in C_i$
if and only if $\alpha$ is a positive root,
$\xi*{[\alpha]}=0$, and $p(\alpha) > p(\beta)+p(\gamma)+\dots$ for
any nontrivial decomposition of $\alpha$ as a sum of positive
roots $\alpha = \beta+\gamma+\dots$ with $\xi*{[\beta]} =
\xi*{[\gamma]} = \dots = 0$.
\end{theorem}

What about the multiplicative equation?
By the Riemann-Hilbert correspondence, any solution
arises as the monodromy of the
differential equation given by a logarithmic connection
on a vector bundle for $\PP^1$. (Note that a Fuchsian ODE (\ref{e:ode}),
as hoped for in Hilbert's 21st problem
in his 1900 ICM talk, will not suffice,
nor will a Fuchsian system of linear differential equations,
or equivalently  a logarithmic connection on a trivial vector bundle,
as discussed by Bolibruch in his ICM talk~\cite{Bolibruchicm}.)
Now a theorem of Weil \cite{Weil} asserts that a
vector bundle on a compact Riemann surface has a holomorphic
connection if and only if its indecomposable direct summands
have degree 0. There is an analogous theorem for parabolic
bundles and compatible logarithmic connections, see \cite{CBipb},
and, using it, Theorem~\ref{t:kaccoh} implies the following.

\begin{theorem}[\cite{CBipb}]
\label{t:conjclo}
There is a solution to $A_1\dots A_k = 1$ with $A_i\in \overline{C_i}$
if and only if
$\alpha$ can be written as a sum of positive roots $\alpha = \beta+\gamma+\dots$ with
$\xi^{[\beta]} = \xi^{[\gamma]} = \dots = 1$.
\end{theorem}

The \emph{Deligne-Simpson problem}, see \cite{Kostovsurvey},
asks when there is an irreducible solution to (\ref{e:prod}) with $A_i\in C_i$.
By considering multiplicative analogues
of preprojective algebras, Crawley-Boevey and Shaw
deduce the following from Theorem~\ref{t:conjclo}.

\begin{theorem}[Crawley-Boevey and Shaw \cite{CBSh}]
\label{t:dsp}
For there to be an irreducible solution to
$A_1\dots A_k = 1$ with $A_i\in C_i$
it is sufficient that $\alpha$ be a positive root,
$\xi^{[\alpha]}=1$, and $p(\alpha) > p(\beta)+p(\gamma)+\dots$ for
any nontrivial decomposition of $\alpha$ as a sum of positive
roots $\alpha = \beta+\gamma+\dots$ with $\xi^{[\beta]} =
\xi^{[\gamma]} = \dots = 1$.
\end{theorem}

The condition in the theorem has now also been
shown to be necessary \cite{CBds}.
For some recent work related to multiplicative preprojective algebras,
see \cite{EOR} and~\cite{VdBdp}.

\end{document}